\documentclass[12pt]{amsart}  
\usepackage{amsmath,amsthm}  
\newcommand{\Z}{\mathbb{Z}}  
\newcommand{\N}{\mathbb{N}}  
\newcommand{\Pb}{{\bf P}}  
\newcommand{\Rb}{{\bf R}}

\newcommand{\reg}{{\rm reg}}

\newcommand{\projdim}{{\rm proj\,dim}}

\newcommand{\hdeg}{{\rm hdeg}}  
\newcommand{\Image}{{\rm Im}}  
\newcommand{\Ker}{{\rm Ker}}  
  
\newcommand{\Hom}{{\rm Hom}}  
\newcommand{\Tor}{{\rm Tor}}

\newcommand{\supp}{{\rm supp}}

\newcommand{\dirsum}{{\oplus}}  
\newcommand{\union}{{\cup}}  
  
\newcommand{\tensor}{{\otimes}}

\newcommand{\pnt}{{\raise 0.5mm \hbox{\large\bf.}}}  
\newcommand{\lpnt}{{\hbox{\large\bf.}}}  
\newcommand{\polyseq}{{poly-exact-sequence}}  
\newcommand{\Sym}{{\rm Sym}}

\newtheorem{Theorem}{\bf Theorem}[section]  
\newtheorem{Lemma}[Theorem]{\bf Lemma}  
\newtheorem{Corollary}[Theorem]{\bf Corollary}  
\newtheorem{Proposition}[Theorem]{\bf Proposition}  
\newtheorem{Remark}[Theorem]{\bf Remark}  
  
\newtheorem{Example}[Theorem]{\bf Example}  
\newtheorem{Example and Definition}[Theorem]{\bf Example and Definition}  
  
\newtheorem{Definition}[Theorem]{\bf Definition}

\title{ Betti numbers of $\mathbb{Z}^n$-graded modules}  
  
\author{Morten Brun} 
\address{FB Mathematik/Informatik, Universit\"at Osnabr\"uck, 49069 Osnabr\"uck, Germany}
\email{brun@mathematik.uni-osnabrueck.de}

\author{Tim R\"omer} 
\address{FB Mathematik/Informatik, Universit\"at Osnabr\"uck, 49069 Osnabr\"uck, Germany}
\email{troemer@mathematik.uni-osnabrueck.de}

\begin{document}  
  
\begin{abstract}  
Let $S=K[X_1,\ldots,X_n]$ be the polynomial ring over a field $K$.  
For bounded below $\Z^n$-graded $S$-modules $M$ and $N$ we show that  
if $   
\Tor^S_p(M,N) \neq 0$, then for $0 \le i \le p$, the dimension of the  
$K$-vector space    
$\Tor^S_i(M,N)$ is at least $\binom p i$.  
In particular,  
we get lower bounds for the total Betti numbers of such modules. These results are   
related to a conjecture of Buchsbaum and Eisenbud.
\end{abstract}  
\maketitle  
%
%
%
\section{Introduction}  
Let $S=K[X_1,\ldots,X_n]$ be the polynomial ring over a field $K$  
with $n$-variables.  
Observe that $S$ has both a natural $\mathbb{Z}$-grading and a natural  
$\mathbb{Z}^n$-grading  
by setting $\deg(X_i)=1$ or $\deg(X_i)=\varepsilon_i$ respectively  
where $\varepsilon_i$ is $i$-th canonical basis vector of $\mathbb{Z}^n$.  
Given $a = (a_1,\dots,a_n) \in \Z^n$ we let $|a| = a_1 + \dots +a_n$.  
We call a $\Z^n$-graded $S$-module $M$ {\em bounded below} if there  
exists $c \in \Z$ such   
that $M_a = 0$ for $|a| \le c$.  In particular, every 
finitely generated 
$\Z^n$-graded $S$-module is bounded below. 
Our main result is a strengthening of a result of Charalambous
concerning $S$-modules of finite length \cite[Theorem 2]{CH91}:
\begin{Theorem}  
\label{mainthm0}  
Let $M, N$ be bounded below $\Z^n$-graded $S$-modules.  
If $\Tor^S_p(M,N) \ne 0$,  
then for $0 \le i \le p$ the $K$-vector space dimension of
$\Tor^S_i(M,N)$ is at least $\binom p i$.   
\end{Theorem}
  
This result is related to a conjecture of Buchsbaum and Eisenbud.
Given an $S$-module $M$ and $i \in \mathbb{N}$ we denote with
$\beta^S_{i}(M)=\dim_K \Tor^S_i(M,K)$ the $i$-th Betti number of $M$.
In this situation Buchsbaum and Eisenbud conjectured the following:
\vspace{3pt}  
\\  
{  
\label{buchsbaum_eisenbud}   
{\bf Conjecture} (Buchsbaum-Eisenbud \cite{BUEI77}):  
\it  
If $M$ is an $S$-module of finite length, then  
$$  
\beta^S_{i}(M)\geq \binom{n}{i} \text{ for } i=0,\ldots,n.  
$$
}
\vspace{3pt}  

For a finite length module $M$ over a regular local ring 
the conjecture was submitted to Hartshorne's problem list by G.\ Horrocks \cite{HA79}.
For this reason, the conjecture is often referred to as Horrocks' Problem. 
However, we will refer to it as the Buchsbaum-Eisenbud conjecture,
since their article appeared before Hartshorne's problem list.
The lower bounds for the Betti numbers imply in particular that  
$\sum_{i=0}^n \beta^R_{i}(M) \geq 2^n$, which was studied for example   
in \cite{AVBU93}.  

In a series of papers (\cite{CH91}, \cite{EVGR88},  \cite{SA90}) the  
Buchsbaum-Eisenbud conjecture was verified   
in the case where $M$ is a $\mathbb{Z}^n$-graded $S$-module of finite length.   

Given a finitely generated $\mathbb{Z}$-graded $S$-module $M$ 
we denote  
with $\projdim(M)=\max\{i\in \mathbb{N}\colon \beta_i^S(M) \neq 0\}$  
the {\em projective dimension} of $M$. 
The module $M$ is said to have a {\em linear} resolution if
$M$ has a minimal $\mathbb{Z}$-graded free resolution
with homogeneous homomorphisms (of degree $0$):
\begin{displaymath}
  0 \to S(-p)^{\beta_p} \to
  S(-(p-1))^{\beta_{p-1}}\to \dots \to S(-1)^{\beta_1} 
  \to S^{\beta_0}  \to M \to 0.
\end{displaymath}
The following result of Herzog
and K\"uhl suggests that the statement of the Buchs\-baum-Eisenbud
conjecture may be extended to modules that are not necessarily of finite
length: 
\begin{Theorem}[\cite{HEKU84}]   
\label{heku}
If a finitely generated $\Z$-graded $S$-module $M$ has a
linear resolution, then
$$  
\beta^S_{i}(M)\geq \binom{\projdim(M)}{i} \text{ for } i=0,\ldots,\projdim(M).  
$$   
\end{Theorem}

However, as a result of Bruns shows \cite[Korollar 1 zu Satz 2]{Bruns},
the above lower bounds will not hold without restrictions on $M$.
In \cite[Theorem 3]{CH91} the 
estimate $\beta^S_{i}(M)\geq \binom{n-d}{i}$ is proved for a
finitely generated $\Z^n$-graded $S$-module $M$, where $d$
is the dimension of $M$. 
Our main result \ref{mainthm0} implies the sharper bounds 
of \ref{heku} in the case that $M$ is a
bounded below $\Z^n$-graded $S$-module. 

The paper is organized as follows. In Section \ref{section1} we define  
the concept of a \polyseq. A \polyseq \ is a 
collection of exact sequences well-suited for inductive construction  
of elements. This inductive procedure is formulated in Theorem
\ref{main_thm_poly}.
In Section \ref{section2} we study a particular  
\polyseq \ related to Theorem \ref{mainthm0}. Theorem \ref{mainthm0}
then follows as a special case of Theorem \ref{main_thm_poly}. In Section
\ref{betti-bounds} we specialize our discussion to Betti numbers.
%
%
%
\section{\polyseq s}  
\label{main_section}  
\label{section1}  
In this section we introduce our concept of a \polyseq \ and we study  
their main properties.   
For $a=(a_1,\ldots,a_n) \in \mathbb{Z}^n$   
set $\supp(a)=\{i \colon a_i \neq 0\} \subseteq [n]=\{1,\ldots,n\}$  
and $|a|=\sum_{i=1}^n a_i$.   
We define $b \preceq a$ for $b=(b_1,\ldots,b_n) \in \mathbb{Z}^n$  
if $b_i \leq a_i$ for $i=1,\ldots,n$. Let $\Pb([n])$ denote the power  
set of $[n]$.   
\begin{Definition}\it  
\label{polyseq}  
Let $\bf{Ab}$ denote the category of abelian groups.  
A map  
$$  
H \colon \Pb([n]) \times \mathbb{N} \times \mathbb{Z}^n   
\to  
\bf{Ab}\it, \quad  
(F,i,a) \mapsto H(F,i,a)=H_i(F)_a  
$$  
is called a \it \polyseq\rm, if the following is satisfied:  
\begin{enumerate}  
\item $H_i(\emptyset) = 0$ for $i >0$.  
\item  
$H_i(F)_a=0$ for $|a|\ll 0$.  
\item  
For all $F \subseteq [n]$, $i \in \mathbb{N}$, $a \in \mathbb{Z}^n$  
and $s \in [n]\setminus F$  
there exist group homomorphisms  
\begin{eqnarray*} 
&&\iota(F,i,a,s) \colon H_i(F)_a \to H_i(F \union \{s\})_a, \\ 
&&\pi(F,i,a,s) \colon H_i(F\union \{s\})_a \to H_{i-1}(F)_{a-\varepsilon_s}, \\  
&&\partial(F,i,a,s) \colon H_i(F)_{a-\varepsilon_s} \to H_i(F)_a,  
\end{eqnarray*} 
such that the following sequences are exact:  
$$  
H_i(F)_{a-\varepsilon_s}  
\xrightarrow{\partial_s}  
H_i(F)_{a}  
\xrightarrow{\iota_s}  
H_i(F\union \{s\})_{a}  
\xrightarrow{\pi_s}  
H_{i-1}(F)_{a-\varepsilon_s}  
\xrightarrow{\partial_s}  
H_{i-1}(F)_{a}.  
$$  
To simplify the notation we write $\iota_s$, $\pi_s$ and $\partial_s$   
if it is clear from the context which homomorphism we mean. We refer  
to the above exact sequence as the {\em fundamental exact sequence for $H$ at  
$(F,i,a,s)$}.   
\end{enumerate}  
We set $H_i(F)=\dirsum_{a \in \mathbb{Z}^n}H_i(F)_a$.  
For $y \in H_i(F)_a$ we call $i$ the homological degree ($\hdeg$)  
and $a$ the degree ($\deg$) of $y$.  
\end{Definition}  
  
Note that we use the notation $H_i(F)_a$ for the values of $H$,  
which reminds of homology groups. In fact in our applications   
$H_i(F)_a$ will be the $\mathbb{Z}^n$-degree $a$ part of the  
$i$-th homology group of a $\mathbb{Z}^n$-graded complex. The above  
definition was motivated by the following example:  
\begin{Example}\rm  
Let $S = K[X_1,\dots,X_n]$ be the polynomial ring over a field $K$ and 
let $M$ be   
a finitely generated $\mathbb{Z}^n$-graded $S$-module.   
For  $F \subseteq[n]$ and $i \in \mathbb{N}$  
let $H_i(F,M)$ be the 
$i$-th Koszul homology of $M$ with respect to the  
variables $X_i$ with $i \in F$. It is well-known that   
for $F \subseteq [n]$ and $s \in [n]\setminus F$ there exists  
the long exact Koszul homology sequence:
$$  
\ldots  
\to  
H_i(F,M)(-\varepsilon_s)  
\xrightarrow{\partial_s}  
H_i(F,M)  
\xrightarrow{\iota_s}  
H_i(F \union \{s\},M)  
\xrightarrow{\pi_s}  
$$  
$$  
H_{i-1}(F,M)(-\varepsilon_s)  
\xrightarrow{\partial_s}  
H_{i-1}(F,M)  
\to  
\ldots  
$$ 
where the maps are $\mathbb{Z}^n$-homogeneous (see \cite[section 1.6]{BRHE98} for details).  
The Koszul homology groups are $\mathbb{Z}^n$-graded.  
Thus we may define $H_i(F)_a=H_i(F,M)_a$.  
Now one verifies that all conditions of \ref{polyseq} are satisfied.  
This is a special case of the \polyseq\ studied in \ref{crucial}.  
\end{Example}  
  
In the rest of this section we present    
a technique to give lower bounds for the number of elements in the groups  
$H_i(F)$. Below $H$ denotes a \polyseq.  
\begin{Lemma}  
\label{vanishing}  
If $i > |F|$, then $H_i(F)_a = 0$ for every $a \in \Z^n$.  
\end{Lemma}  
\begin{proof}  
We prove the lemma by induction on $|F|$. By the axiom $(i)$ for a \polyseq \  
the lemma holds in the case where $|F| = 0$. Suppose that the lemma  
holds for $F$. It follows from the fundamental exact sequence for $H$ at  
$(F,i,a,s)$ that the lemma holds for $F \union \{ s \}$.  
\end{proof}  
The following lemma observes that it is possible to 'push' elements at  
a fixed homological degree.  
\begin{Lemma}\it  
\label{crucial1}  
Let   
$F \subseteq[n]$,   
$s \in [n]\setminus F$ and   
$0 \neq y \in H_i(F)_{a}$ for some $0 \leq i \leq |F|$,  $a \in \mathbb{Z}^n$.  
Then there exists an element   
$$  
0 \neq y_s \in H_{i}(F \union \{s\})_{a - b_s \varepsilon_s}  
\text{ for some }   
b_s \in \mathbb{N}.  
$$  
\end{Lemma}  
\begin{proof}  
Define $l = \min\{|c| \colon H_{i}(F)_{c} \neq 0, c \in \mathbb{Z}^n \}$.  
Observe that $l$ is well defined because axiom $(ii)$ ensures that  
$H_{i}(F)_{c}=0$ for $|c|\ll 0$.   
We prove this lemma by induction on $|a|- l$.  
Let $|a|=l$. We consider the fundamental exact sequence for $H$ at  
$(F,i,a,s)$: 
$$  
H_i(F)_{a-\varepsilon_s}  
\xrightarrow{\partial_s}  
H_i(F)_{a}  
\xrightarrow{\iota_s}  
H_i(F\union \{s\})_{a}  
\xrightarrow{\pi_s}  
H_{i-1}(F)_{a-\varepsilon_s}  
\xrightarrow{\partial_s}  
H_{i-1}(F)_{a}.  
$$  
By the definition of $l$ we have $H_{i}(F)_{a - \varepsilon_s}=0$.  
Therefore $\iota_s : H_{i}(F)_{a}  
\to H_{i}(F \union \{s\})_{a}$ is injective and $y_s=\iota_s(y)$  
fulfils the requirements.   
  
Next assume that $|a|>l$.   
We have to consider two cases. If $\iota_s(y)\neq 0$, then we choose   
$y_s=\iota_s(y)$ and $b_s=0$. Otherwise $\iota_s(y)= 0$.   
But by the fundamental exact sequence for $H$ at $(F,i,a,s)$  
there exists a   
$0 \neq y' \in H_{i}(F)_{a - \varepsilon_s}$ with $\partial_s(y')=y$.  
By the induction hypothesis we get an element   
$$  
0 \neq y_s \in H_{i}(F\union \{s\})_{a - \varepsilon_s - b'_s \varepsilon_s}  
$$  
for some $b'_s \in \mathbb{N}$. With $b_s=b'_s+1$ we found again the  
desired element in degree $a - b_s \varepsilon_s$.  
\end{proof}  
  
We are ready to prove the main result of this section.  
Let $i,j \in \mathbb{Z}$. 
Recall the convention that  
$\binom{i}{0}=1$ and that    
$\binom{i}{j}=0$ unless $0 \le j \le i$.
 
\begin{Theorem}\it  
\label{main_thm_poly}  
Suppose given    
a subset $F$ of $[n]$,   
numbers $r$ and $p$ satisfying that $0 \leq r \leq p \leq |F|$, a subset   
$\{l_1,\ldots, l_r\} \subseteq F$ of cardinality $r$ 
and a nonzero element  
$y \in H_p(F)_{a}$ for some $a \in \mathbb{Z}^n$. 
Then for all $i \in \mathbb{N}$ there exist   
$\binom{p-r}{i}$ many $I \subseteq F \setminus \{l_1,\ldots, l_r\}$ and  
nonzero elements $y_I \in H_{p-i}(F)_{a-b_I}$ with  
\begin{enumerate}  
\item  
$b_I \in \N^n \subseteq \Z^n$ and $|b_I| \ge i$.  
\item  
$\supp(b_I)\subseteq F \setminus \{l_1,\ldots, l_r\}$.  
\item  
The degrees $\deg(y_I)$ are pairwise distinct.  
\end{enumerate}  
\end{Theorem}  
\begin{proof}  
The theorem clearly holds in the case where $|F| = r$. In particular it 
holds for $|F| = 0$. Let $0 \le r < |F|$, and assume by induction that 
the theorem holds for every $G \subseteq [n]$ with $|G| < |F|$. In order 
to prove the induction step we consider two cases. 
  
Case 1:  
Assume that for some $s \in F \setminus \{l_1,\ldots, l_r\}$ we have $\pi_s(y)=0$.  
By the fundamental exact sequence for $H$ at $(F\setminus \{ s\},i,a,s)$   
$$  
H_i(F\setminus \{ s\})_{a-\varepsilon_s}  
\xrightarrow{\partial_s}  
H_i(F\setminus \{ s\})_{a}  
\xrightarrow{\iota_s}  
H_i(F )_{a}  
\xrightarrow{\pi_s}  
H_{i-1}(F\setminus \{ s\})_{a-\varepsilon_s}  
\xrightarrow{\partial_s}  
H_{i-1}(F\setminus \{ s\})_{a},  
$$  
there exists an element $0 \neq y' \in H_{p}(F\setminus\{s\})_{a}$  
with $\iota_s(y')=y$.  
Lemma \ref{vanishing} ensures that $p \le |F \setminus \{s\}|$ and  
therefore we can apply the induction hypothesis to  
$y'$ and get     
that for all $i \geq 0$ there exist   
$\binom{p-r}{i}$ many 
$I' \subseteq (F\setminus \{s\}) \setminus \{l_1,\ldots, l_r\}$ and  
$0 \neq y'_{I'} \in H_{p-i}(F\setminus \{s\})_{a-b'_{I'}}$ with  
\begin{enumerate}  
\item  
$b'_{I'} \in \N^n \subseteq \Z^n$, $|b'_{I'}| \ge i$.  
\item  
$\supp(b'_{I'})\subseteq (F\setminus \{s\}) \setminus \{l_1,\ldots, l_r\}$.  
\item  
The degrees $\deg(y'_{I'})$ are pairwise distinct.  
\end{enumerate}  
Let $I=I' \union \{s\}$. By \ref{crucial1} we can push the elements 
$y'_{I'}$ to   
$0 \neq y_I \in H_{p-i}(F)_{a-b'_{I'}- b_s \varepsilon_s}$  
for some $b_s \geq 0$. With   
$b_I=b'_{I'}+b_s \varepsilon_s$ we are done in this case.  
  
Case 2:  
Assume that for all $s \in F \setminus \{l_1,\ldots, l_r\}$ we have $\pi_s(y)\neq 0$.  
  
For $i=0$ define $I=\emptyset$ and $y_I=y$. Let $i \geq 1$.  
Since   
$$  
|F \setminus \{l_1,\ldots, l_r \}|\geq p-r  
$$   
we can choose   
$s_1,\ldots,s_{p-r} \in F \setminus \{l_1,\ldots, l_r \}$.  
By the induction hypothesis the theorem applies to $F\setminus\{s_j\} 
\subseteq [n]$, $0 \le r+j-1 \le p-1 \le |F \setminus \{s_j \}|$, 
$\{l_1,\dots,l_r,s_1,\dots,s_{j-1} \}$ $\subseteq F\setminus\{s_j\}$ and 
$\pi_{s_j}(y) \in H_{p-1}(F\setminus\{s_j\})_{a - 
\varepsilon_{s_j}}$. Hence for each $s_j$ we get 
$\binom{p-1-r-j+1}{i-1}$ many  
$I_{s_j} \subseteq (F \setminus \{s_j\}) \setminus  
\{l_1,\ldots, l_r,s_1,\ldots,s_{j-1}\}$ 
with $|I_{s_j}|=i-1$   
and elements   
$$  
0 \neq y'_{I_{s_j}} \in   
H_{p-1-(i-1)}(F \setminus \{s_j\})_{a-\varepsilon_{s_j}-b'_{I_{s_j}}}$$   
with  
\begin{enumerate}  
\item 
$b'_{I_{s_j}} \in \N^n \subseteq \Z^n$ and $|b'_{I_{s_j}}| \ge i-1$.  
\item  
$\supp(b'_{I_{s_j}})\subseteq (F \setminus \{s_j\}) \setminus 
\{l_1,\ldots, l_r,s_1,\ldots,s_{j-1}\}$.   
\item  
The degrees $\deg(y'_{I_{s_j}})$ are pairwise distinct.  
\end{enumerate}  
By \ref{crucial1} we can push the elements $y'_{I_{s_j}}$ to   
$$  
0 \neq y_{I_{s_j}} \in    
H_{p-i}(F)_{a-(1+b_{s_j})\varepsilon_{s_j} - b'_{I_{s_j}}}  
$$  
for some $b_{s_j}\geq 0$.  
Define $b_{I_{s_j}}=(1+b_{s_j})\varepsilon_{s_j} + b'_{I_{s_j}}$.  
Note that $|b_{I_{s_j}}| \ge i$ and that for $s_j\neq s_{j'}$ with
$j<j'$ we have $\deg(y_{I_{s_j}}) \neq \deg(y_{I_{s_{j'}}})$,   
because   
$$  
s_j \in \supp(b_{I_{s_j}})\subseteq F \setminus \{l_1,\ldots, l_r,s_1\ldots,s_{j-1}\}  
$$  
and  
$$  
s_j \not\in   
\supp(b_{I_{s_j'}})\subseteq F \setminus \{l_1,\ldots, l_r,s_1,\ldots,s_j,\ldots,s_{j'-1}\}.  
$$  
Clearly $\deg(y_{I_{s_j}}) \neq \deg(y_{I'_{s_{j}}})$ for sets $I_{s_j}\neq I'_{s_j}$  
with respect to the same $s_j$.  
We have produced  
$$  
\binom{p-1-r}{i-1}+\binom{p-1-r-1}{i-1}+\cdots+\binom{0}{i-1}=\binom{p-r}{i}  
$$  
many elements $0 \neq y_{I_{s_j}}\in H_{p-i}(F)$   
with pairwise distinct degrees. Thus we are done.  
  
\end{proof}  
  
\section{A \polyseq \ for Tor-groups}
\label{section2}
In this section we give an example of a \polyseq\ and present an
application of   
\ref{main_thm_poly}.   
Let $S=R[X_1,\ldots,X_n]$ be the polynomial algebra over a commutative  
and unital ring $R$  
with $n$-variables equipped with the grading $\deg(X_i)=\varepsilon_i$.   
For $F \subseteq [n]$   
we define $R[F] \subseteq S$ to be the subalgebra generated   
by the variables $X_i$ with $i \in F$.  
  
The following is a lemma in \cite{RO79} (Section 9).   
Since it is a crucial observation for our paper we reproduce its  
proof.  
\begin{Lemma}\it  
\label{lem_1}  
Let $M$ be a $\mathbb{Z}^n$-graded $S$-module,   
$F \subseteq [n]$ and $s \in [n]\setminus F$. Then  
there exists a short exact sequence of $\mathbb{Z}^n$-graded $S$-modules 
of the form 
$$  
0  
\to  
S \tensor_{R[F]} M (-\varepsilon_s)  
\xrightarrow{\partial_s}  
S \tensor_{R[F]} M   
\xrightarrow{\iota_s}  
S \tensor_{R[F \union \{s\}]} M  
\to  
0  
$$  
with  
$\partial_s(1 \tensor m)=X_s \tensor m  - 1 \tensor X_s m$  
and  
$  
\iota_s$ induced by the inclusion $R[F] \subseteq R[F \union \{s\}]$.
\end{Lemma}  
\begin{proof}  
Since $M$ is a $\mathbb{Z}^n$-graded $S$-module,  
it is a $\mathbb{Z}^n$-graded $R[F \union \{s\}]= R[F][X_s]$-module.  
We claim that the following sequence is exact   
\begin{equation}  
\label{ex_ro_seq}  
0  
\to  
R[F \union \{s\}] \tensor_{R[F]} M (-\varepsilon_s)  
\xrightarrow{\partial_s}  
R[F \union \{s\}] \tensor_{R[F]} M   
\xrightarrow{\iota_s}  
M  
\to  
0  
\end{equation}  
with $\partial_s(1 \tensor m)=X_s \tensor m  - 1 \tensor X_s m$  
and $\iota_s(r \tensor m)=rm$.  
  
Clearly $\iota_s$ is surjective and $\iota_s \circ \partial_s=0$.  
We claim that $\Ker(\iota_s)=\Image(\partial_s)$.   
If $\sum_{i=0}^k X_s^i \tensor m_i \in \Ker(\iota_s)$, then  
$  
0=\sum_{i=0}^k X_s^i m_i.  
$  
Define   
$$  
n_j=\sum_{i=j+1}^k X_s^{i-j-1} m_i \text{ for } j=0,\ldots, k-1.  
$$  
Then we have $m_k=n_{k-1}$, $m_j=n_{j-1}-X_sn_j$ for $j=1,\ldots,k-1$ 
and $m_0=-X_sn_0$.   
Hence  
$  
\partial_s(\sum_{j=0}^{k-1} X_s^j \tensor n_j)=\sum_{i=0}^k X_s^i \tensor m_i.  
$  
It remains to show that $\partial_s$ is injective.   
Assume that $\sum_{i=0}^k X_s^i \tensor m_i \in \Ker(\partial_s)$.  
It follows from  
$$  
0  
=  
\partial(\sum_{i=0}^k X_s^i \tensor m_i)  
=  
\sum_{i=0}^k (X_s^{i+1} \tensor m_i - X_s^i \tensor X_sm_i)  
$$  
that  
$  
0=m_k=m_{k-1}-X_sm_k=\ldots= m_{0}-X_s m_1  
$  
and therefore $m_i=0$ for $i=0,\ldots,k$.  
  
Since $S$ is a flat $R[F \union \{s\}]$-module   
we can apply $S \tensor_{R[F\union \{s\}]} $ to the sequence (\ref{ex_ro_seq})  
and obtain the exact sequence  
$$  
0  
\to  
S \tensor_{R[F]} M (-\varepsilon_s)  
\xrightarrow{\partial_s}  
S \tensor_{R[F]} M   
\xrightarrow{\iota_s}  
S \tensor_{R[F \union \{s\}]} M  
\to  
0  
$$  
of $\mathbb{Z}^n$-graded $S$-modules.  
\end{proof}  
  
In the proof of \ref{lem_1} we did not use the structure  
of $\mathbb{Z}^n$-graded modules. This result holds for arbitrary $S$-modules.  
  
Recall that we say that a $\mathbb{Z}^n$-graded $S$-module is \it  
bounded below \rm if $M_a=0$   
for $a \in \mathbb{Z}^n$ and $|a| \ll 0$.   
  
\begin{Proposition}\it  
\label{crucial}  
Let $M, N$ be $\mathbb{Z}^n$-graded $S$-modules.   
If $M$ and $N$ are bounded below and either $M$ or $N$ is flat  
considered as an $R$-module,  
then $T$     
with  
$$  
T_i(F)_a = \Tor_i^S(S \tensor_{R[F]} M, N)_a  
$$  
for $F \subseteq [n], i \in \mathbb{N}$ and $a \in \mathbb{Z}^n$ is a \polyseq.  
\end{Proposition}  
\begin{proof}  
We have to verify the conditions of a \polyseq\ in \ref{polyseq}.  
  
(i): Since either $M$ or $N$ is flat considered as an $R$-module we have  
$\Tor^R_i(M,N) = 0$ for $i \ne 0$. The axiom now follows from the flat  
base-change formula: $T_i(\emptyset)_a= \Tor^S_i(S \otimes_R M,N) \cong  
\Tor^R_i(M,N)$. (See for example Weibel \cite{WE94}, Proposition 3.2.9.)  
  
(ii): $T_i(F)$ is a $\mathbb{Z}^n$-graded $S$-module, which is bounded below,  
because $S \tensor_{R[F]} M$ and  $N$ are bounded below. Thus  
$T_i(F)_a=0$ for $|a|\ll 0$.   
  
(iii):  
By \ref{lem_1} we have   
for  all $F \subseteq[n]$, $i \in \mathbb{N}$, $a \in \mathbb{Z}^n$  
the exact sequence  
$$  
0  
\to  
S \tensor_{R[F]} M (-\varepsilon_s)  
\xrightarrow{\partial_s}  
S \tensor_{R[F]} M   
\xrightarrow{\iota_s}  
S \tensor_{R[F \union \{s\}]} M  
\to  
0.  
$$  
The degree $a$-part of the corresponding $\mathbb{Z}^n$-graded long exact Tor sequences   
is  
$$  
\cdots  
T_i(F)_{a-\varepsilon_s}  
\xrightarrow{\partial_s}  
T_i(F)_{a}  
\xrightarrow{\iota_s}  
T_i(F \union \{s\})_{a}  
\xrightarrow{\pi_s}  
T_{i-1}(F)_{a-\varepsilon_s}  
\xrightarrow{\partial_s}  
T_{i-1}(F)_{a}  
\cdots  
$$  
Here $\pi_s$ is the connecting homomorphism and the other maps are induced by  
the corresponding ones in \ref{lem_1}. This concludes the proof.  
\end{proof}  
  
The following theorem implies Theorem \ref{mainthm0} and  
\cite[Theorem 2]{CH91}.   
\begin{Theorem}\it  
\label{buei_before}  
Let $M, N$ be bounded below $\mathbb{Z}^n$-graded $S$-modules and assume  
that either $M$ or $N$ is flat considered as an $R$-module.   
If $\Tor_p^S(M,N)_a \neq 0$ for some $a \in \mathbb{Z}^n$, then for $i  
\in \{ 0, \dots, p \}$ there  
exist at least $\binom p i$ distinct $b \in \Z^n$ with $b \preceq a$,
$|b| \le |a| -p + i$  
and $\Tor_i^S(M,N)_b \ne 0$.  
In particular, if $R=K$ is a field, we have $\dim_K \Tor_i^S(M,N) \geq  
\binom{p}{i}$.   
\end{Theorem}  
\begin{proof}  
By \ref{crucial} we get the \polyseq\ $T$.   
Note that $T_i([n])=\Tor_i^S(M,N)$.  
Since $T_p([n])_a \neq 0$ we get an element $0 \neq y \in T_p([n])_a$.  
The result now follows directly from \ref{main_thm_poly} with $F =  
[n]$ and $r = 0$.   
\end{proof}  

\begin{Remark}\rm
  In the above theorem it is essential to require that $M$ and $N$ are
bounded below. In order to see this, let $R=K$ be a field and $n=1$. 
Then $S = K[X]$ is the polynomial ring with one indeterminate. 
Since $\Tor^{S}_i(K[X,X^{-1}],K) = 0$ for every
$i \geq 0$ the $S$-module $M = K[X,X^{-1}]/S X$ has 
$\Tor_1^S(M,K) \cong K$ and $\Tor^S_0(M,K)=0$.
\end{Remark}
\section{Lower bounds for Betti numbers}  
\label{betti-bounds}    
In this section we let $S = K[X_1,\dots,X_n]$ be the
polynomial ring over a field $K$ with $n$-variables.  
Let $M$ be a $\mathbb{Z}^n$-graded $S$-module.  
For $i \in \mathbb{N}$ and $a \in \mathbb{Z}^n$ we denote with  
$\beta_{i,a}^S(M)=\dim_K \Tor_i(M,K)_a$ the $\mathbb{Z}^n$-graded  
Betti numbers of $M$   
and   
with  
$\beta_{i}^S(M)=\sum_{a \in \mathbb{Z}^n} \beta_{i,a}^S(M)$ the total Betti numbers of $M$.  
  
The following result verifies the bounds of \ref{heku}
in the case of   
$\mathbb{Z}^n$-graded modules.   
\begin{Corollary} 
\label{betti_cor} 
\it  
Let $M$ be a bounded below $\mathbb{Z}^n$-graded $S$-module.  
If $\beta^S_{p,a}(M) \neq 0$ for some $0 \leq p \leq n$ and $a \in \mathbb{Z}^n$, then  
$$  
\sum_{b \in \mathbb{Z}^n,\ b \preceq a,\ 
|b| \le |a| - p + i} \beta^S_{i,b}(M) \geq \binom{p}{i} \text{ for } i=0,\ldots,p.  
$$  
In particular, we have:  
\begin{enumerate}  
\item  
$\beta^S_{i}(M) \geq \binom{p}{i}$.  
\item  
$\sum_{i=0}^p \beta^S_i(M) \geq 2^p$.  
\end{enumerate}  
\end{Corollary}  
\begin{proof}  
Apply \ref{buei_before} with $N=K$.  
\end{proof}  
  
Let $M$ be a finitely generated $\mathbb{Z}^n$-graded $S$-module.  
For $i,j \in \mathbb{Z}$   
we denote with $\beta_{i,j}^S(M)=\sum_{a \in \mathbb{Z}^n, |a|=j}\beta_{i,a}^S(M)$  
the $\Z$-graded Betti numbers of $M$.  
  
The Castelnuovo-Mumford regularity for a finitely generated $\mathbb{Z}$-graded $S$-module $0 \neq M$   
is defined as  
$\reg(M)=\max\{j \in \mathbb{Z}\colon \beta^S_{i,i+j}(M)\neq 0 \text{ for some } i \in \mathbb{N} \}.$  
For $k \in \{0,\ldots,n\}$ we  
define $d_k(M)=\min(\{j \in \mathbb{Z} \colon \beta^S_{k,k+j}(M)\neq 0 \} \union \{\reg(M)\})$.   
Let $(F_\lpnt,d_\lpnt)$ 
be the minimal $\Z$-graded free resolution of $M$. Since
$F_i = \bigoplus_{j \in \Z} S(-i-j)^{\beta^S_{i,i+j}(M)}$ and $d_i(F_i)
\subseteq (X_1,\dots,X_n) F_{i-1}$ we have that 
$d_0(M) \leq d_1(M) \leq \cdots$.  
We are interested in the   
numbers $\beta_i^{k,lin}(M)=\beta^S_{i,i+d_k(M)}(M)$ for  
$i \geq k$. The numbers $\beta_i^{lin}(M)=\beta_i^{0,lin}(M)$ are the Betti numbers of the linear Strand of $M$.  
The second author proved in \cite{RO02} the next result, which we now get easily  
from \ref{betti_cor}.  
  
\begin{Theorem}\it  
\label{syzallg}  
Let $k \in [n]$ and $M$ be a finitely generated $\mathbb{Z}^n$-graded $S$-module.  
If $\beta_p^{k,lin}(M) \neq 0$ for some $p \geq k$, then   
$$  
\beta_i^{k,lin}(M) \geq \binom{p}{i}  
\text{ for } i=k,\ldots, p.  
$$  
\end{Theorem}  
\begin{proof}  
If $\beta_p^{k,lin}(M) \neq 0$, there exists an element  
$0 \neq y \in \Tor_p^S(M,K)_{a}$ with $a \in \mathbb{Z}^n$ and $|a|=p+d_k(M)$.  
By \ref{betti_cor} we have  
$$  
\sum_{b \in \mathbb{Z}^n,\ b \preceq a,\ |b| \le |a| - p + i}  
\beta^S_{i,b} 
\geq \binom{p}{i} \text{ for } i=k,\ldots,p.  
$$  
Since all considered $b$ in this sum have $|b|\leq |a| -p + i=i +d_k(M)$ and   
$\Tor_i^S(M,K)_{c}=0$ for $|c|<i+d_k(M)$ we have  
$\Tor_i^S(M,K)_b \neq 0$ only if $|b|=i+d_k(M)$.  
Hence  
$$  
\beta_i^{k,lin}(M)= 
\sum_{b \in \Z^n,\ |b| = i + d_k(M)} \beta^S_{i,b}(M)  
\geq \binom{p}{i} \text{ for } i=k,\ldots,p.  
$$  
  
\end{proof}  
  

\begin{Corollary}\it  
\label{kth}  
Let $M$ be a finitely generated $\mathbb{Z}^n$-graded $S$-module   
and suppose that $M$ is the $k$-th syzygy module   
in a minimal $\mathbb{Z}^n$-graded free resolution.   
If $\beta_p^{lin}(M) \neq 0$ for some $p \in \mathbb{N}$, then   
$$  
\beta_i^{lin}(M) \geq \binom{p+k}{i+k} \text{ for } i=0,\ldots, p.$$  
\end{Corollary}  
\begin{proof}  
The result follows from Theorem \ref{syzallg} by the fact that  
if $M$ is a $k$-th syzygy module of $N$, then $\beta_i^{lin}(M) =  
\beta_{i+k}^{k,lin}(N)$. For details we refer to    
\cite{RO02} where the result is proved by a similar method.  
\end{proof}  

\begin{Remark}\rm
(a)
The result \ref{kth} verifies the $\Z^n$-graded case of a conjecture of  
Herzog \cite{HE98} on finitely generated $\mathbb{Z}$-graded $S$-modules.  
This conjecture is motivated by a result of Green \cite{GR84}  
(see also Eisenbud and Koh \cite{EIKO91} and Green \cite{GR99})
that contains the case $i=0, k=1$. 
See \cite{RO02} for details and other known cases.

(b) 
Let $W$ be an $n$-dimensional $K$-vector space and $V=\Hom_K(W,K)$ 
the dual vector space. We give the elements of $W$ the degree
1. Thus the elements of $V$ have the degree $-1$. 
Then $S$ can be considered as the symmetric algebra $\Sym(W)$
of $W$. Let $E=\wedge V$ be the exterior algebra of $V$.
A special case (\cite{EIFLSC}, Prop. 2.1) 
of the so-called BGG correspondence shows 
that there exists a functor $\Rb$ 
which is an equivalence between the category 
of graded left $S$-modules  
and the category of linear free complexes over $E$.
(See \cite{EIFLSC} for unexplained terms and proofs.)
Moreover, if $M$ is a graded $S$-module and $\Rb(M)$ is the corresponding
linear free complex over $E$, then we have
$
H^j(\Rb(M))_{i+j}\cong \Tor_i^S(M,K)_{i+j}
$
(\cite{EIFLSC}, Prop. 2.3).
This implies immediately that 
$$
\min\{j \in \Z \colon H^j(\Rb(M))\neq 0\}
=
\min\{j \in \Z \colon \beta^S_{0,j}(M) \neq 0\}
=
d_0(M).
$$
We get that $H^{d_0(M)}(\Rb(M))$ is a graded $E$-module
with $\dim_K H^{d_0(M)}(\Rb(M))_{i+d_0(M)}$ $= \beta^{lin}_{i}(M)$.
Thus Corollary \ref{kth} gives a lower estimate for the Hilbert function
of $H^{d_0(M)}(\Rb(M))$.
\end{Remark}

%
%
%
  
\end{document}